%----------------------------------------------------------------------
%Proposed final version 2/17/97
\magnification 1200
\input amstex
\documentstyle{amsppt}
%\hoffset 5.5truemm
\topmatter
\NoBlackBoxes
\rightheadtext{Local biholomorphisms of real hypersurfaces}
\leftheadtext{M. S. Baouendi,  P. Ebenfelt, and
L. P. Rothschild}
\title Parametrization of local biholomorphisms of real
analytic hypersurfaces\endtitle
\author M. S. Baouendi,  P. Ebenfelt, and
Linda Preiss Rothschild
\endauthor
\address Department of Mathematics, University of California
at San Diego, La Jolla, CA 92093\endaddress
\email sbaouendi\@ucsd.edu, lrothschild\@ucsd.edu\endemail
\address Department of Mathematics, Royal Institute of
Technology, 100 44 Stockholm, Sweden\endaddress
\email ebenfelt\@math.kth.se\endemail
\thanks The first and the third authors are partially
supported by National Science Foundation grant DMS
95-01516. The second author is supported by a grant from
the Swedish Natural Science Research Council.\endthanks
%\abstract \endabstract
%\keywords \endkeywords
%\address Department of Mathematics, University of
%California at
%San Diego,
%La Jolla, CA 92093\endaddress
%\email pebenfel\@euclid.ucsd.edu\endemail
%\address Department of Mathematics, University of
%California at
%San Diego,
%La Jolla, CA 92093\endaddress
%\email sbaouendi\@ucsd.edu\endemail
%\address Department of Mathematics, University of
%California at
%San Diego,
%La Jolla, CA 92093\endaddress
%\email lrothschild\@ucsd.edu\endemail
%\date{\number\year-\number\month-\number\day}\enddate
\loadeufm
\def\tr{\text{\rm tr}}
\def\xo {{x_0}}
\def\sg{\sigma}
\def \emxy{E_{(M,M')}(X,Y)}
\def \semxy{\scrE_{(M,M')}(X,Y)}
\def \jkxy {J^k(X,Y)}
\def \gkxy {G^k(X,Y)}
\def \exy {E(X,Y)}
\def \sexy{\scrE(X,Y)}

\def\prt#1{{\partial \over\partial #1}}
\def\det{{\text{\rm det}}}
\def\wob{{w\over B(z)}}
\def\co{\chi_1}
\def\po{p_0}
\def\fb {\bar f}
\def\gb {\bar g}
\def \qb {\bar Q}
\def \t {\tau}
\def\z{\zeta}
\def \L {\Lambda}
\def\b {\beta}
\def\a {\alpha}
\def\l {\lambda}
%\define \p{\eufm p}
%\define \q{\eufm q}
%\define \m{\eufm m}
%\define \ass{\text{\rm Ass }}
%\define \sol{\text{\rm Sol}}
%\define \atan{\text{\rm arctan }}
\define \im{\text{\rm Im }}
\define \re{\text{\rm Re }}

\define \bR{\Bbb R}
%\define \bR{{\bold R}}
\define \bC{\Bbb C}
%\define \bC{{\bold C}}

%\define \bL{{\bold L}}
\define \scrF{\Cal F}
\define \scrK{\Cal K}
\define \scrM{\Cal M}
\define \scrR{\Cal R}
\define \scrJ{\Cal J}
\define \scrA{\Cal A}

\define \scrL{\Cal L}
\define \scrE{\Cal E}
%\define \ssneq{\subsetneqq}
%\define \ssneq{\subset}
\define \hol{\text{\rm hol}}

\define \Aut{\text{\rm Aut}}

\def\jet#1#2{J^{#1}_{#2}}
\def\gp#1{G^{#1}}
\def\gpo{\gp {2k_0}_0}
\def\emmp {\scrF(M,p;M',p')}
\endtopmatter

\document
\heading \S 0 Introduction\endheading

A real analytic hypersurface $M\subset \bC^N$ is said to be
{\it finitely nondegenerate} at
$p_0\in M$ if there is a positive integer $k$  such that
$$
\text{span }\left\{L^\alpha\rho_Z(p_0,\bar
p_0)\:|\alpha|\leq k\right\}=\bC^N.
\tag0.1
$$ Here $L_1,...,L_{N-1}$ is a  basis for the CR vector
fields on $M$ near
$p_0$, $L^\alpha$ denotes
$L_1^{\alpha_1}...L_{N-1}^{\alpha_{N-1}}$,
$\rho(Z,\bar Z)=0$ is a defining equation of $M$ near
$p_0$, and $\rho_Z=\partial\rho/
\partial Z$. $M$ is called {\it $k$-nondegenerate} at $\po$
if
$k$ is the smallest integer for which (0.1) holds. This
integer is an invariant of
$M$ at
$p_0$; it does not depend on the basis $L_1,...,L_{N-1}$,
the defining function $\rho(Z,\bar Z)$, or the coordinates
used (see [BHR] in which this notion was introduced). The
hypersurface $M$ is $1$-nondegenerate  at
$\po$ if and only if it is Levi nondegenerate at $\po$.

A connected hypersurface $M$ is called {\it holomorphically
nondegenerate} if there is no nontrivial germ of a
holomorphic vector field (i.e. a holomorphic section of
$T^{(1,0)}(\bC^N)$) tangent to $M$. One can show that it
suffices to check this condition at any one point of
$M$.  Furthermore, a connected real analytic hypersurface
$M$ is holomorphically nondegenerate if and only if it is
finitely nondegenerate at some point (and hence at all
points outside a proper real analytic subset of $M$). We
refer the reader to [BHR], [BER1] for details. (See also
Stanton [S2].)

In this paper we study germs of biholomorphisms mapping one
real analytic hypersurface into another and show that at
points of finite nondegeneracy there is a natural parameter
space which has the structure of a real analytic manifold.

Some of the main results of this paper can be described as
follows. Let $M$  and $M'$ be two real analytic
hypersurfaces $k_0
$-nondegenerate at $p$ and $p'$ respectively.  We consider
the (possibly empty) set of  germs
$H$ at $p$ of invertible holomorphic mappings in $\bC^N$
with $H(M) \subset M'$ and $H(p) =p'$. This is a subset of
the space of all germs at $p$ of holomorphic mappings
$\bC^N\to\bC^N$ and as such inherits a natural inductive
limit topology corresponding to uniform convergence on
compact neighborhoods of $p$.   We show (Theorem 1) that
this set, equipped with its natural topology, is mapped
into the complex jet group of order
$2k_0$ of $\bC^N$ homeomorphically  onto a real analytic,
totally real submanifold $\Sigma_{p,p'}$.
Moreover,
$\Sigma_{p,p'}$ is given by equations that can be
explicitly constructed from defining equations of $M$ and
$M'$. In particular, the germs of hypersurfaces $(M,p)$ and
$(M',p')$ are biholomorphically equivalent if and only if
these equations have a solution; we say that $(M,p)$ and
$(M',p')$ are  biholomorphically equivalent if there is a
germ
$H$ of a biholomorphism mapping $M$ into
$M'$ with $H(p) = p'$. We say that $(M,p)$ and
$(M',p')$ are {\it formally equivalent} if the components of
$H$ are assumed only to be formal power series at $p$.   We
show (Theorem 5) that if
$(M,p)$ and
$(M',p')$ are finitely nondegenerate, then they are formally
equivalent if and only if they are  biholomorphically
equivalent.  In fact, we show that every formal equivalence
is a biholomorphic equivalence.

 We denote by $\Aut(M,p)$ the germs of biholomorphisms of
$\bC^N$ at $p$ that fix $p$ and map
$M$ into itself. This is a group under composition,  called
the {\it stability group}  of $M$ at
$p$. It follows from the above results with $M=M'$ and
$p=p'$ that
$\Aut(M,p)$, equipped with its natural topology, is a real
Lie group  provided that $M$ is finitely nondegenerate at
$p$ (Corollary 1.2).  Hence, if
$M$ is holomorphically nondegenerate, then $\Aut(M,p)$ is a
real Lie group for all
$p\in M$ outside a proper real analytic subset of $M$.

If $M$ is Levi nondegenerate, i.e. $1$-nondegenerate, at
$\po$, then the fact that $\Aut (M,\po)$ is a real Lie
group follows from the work of Chern-Moser [CM].  In the
general case, the above result for $\Aut (M,\po)$ is
related to recent work of Stanton [S2] and the present
authors [BER2]. Denote by
$\hol(M,\po)$ the Lie algebra of real  analytic
infinitesimal CR automorphisms of $M$ at $p_0$; i.e.
$\hol(M,\po)$ consists of all germs at $\po$ of vector
fields
$X$ tangent to $M$ such that the local 1-parameter group of
transformations generated by $X$ are biholomorphic
transformations of $\bC^N$ preserving
$M$.  Stanton [S2] proved that
$\hol(M,\po)$ is finite dimensional for every $p_0\in M$ if
and only if $M$ is holomorphically nondegenerate. (A
different proof of this result as well as a proof for the
corresponding result  for merely smooth ($C^\infty$)
infinitesimal CR automorphisms and also higher
codimensional analogs were given in [BER2]).
 Let $\hol_0(M,\po)$ be the elements of
 $\hol(M,\po)$ that vanish at $p_0$. Under
the assumption that $M$ is finitely  nondegenerate at $\po$,
$\Aut(M,\po)$ is a Lie group, as mentioned above, and
$\hol_0(M,\po)$ is its Lie algebra. The finite
dimensionality of
$\hol(M,\po)$ implies that  there is a unique topology on
$\Aut(M,\po)$, considered  as an abstract group, such
that it is a Lie transformation group with $\hol_0(M,\po)$
as its Lie algebra (see e.g. [Ko, p. 13]). The results above
imply that if $M$ is finitely nondegenerate at $\po$ then
this topology must coincide with the natural topology on
$\Aut(M,\po)$.

If
$M$ is a holomorphically nondegenerate hypersurface, then
there is a number
$\ell(M)$, the {\it Levi type} of
$M$,  with $1\leq \ell(M)\leq N-1$, such that $M$ is
$\ell(M)$-nondegenerate at all points outside a real
analytic subset. (See [BHR] and [BER1].)
Theorem 3 deals with the dependence of
$\Sigma_{p,p'}$ on the base points
$p$ and
$p'$ when $M$ and $M'$ are $\ell(M)$-nondegenerate at all
points, and more generally when $M$ and $M'$ are at most
$k_0$-nondegenerate at all points for some $k_0$. We prove
that the union
$\Sigma$ of the submanifolds $\Sigma_{p,p'}$  forms a real
analytic subset of the appropriate jet bundle. In general,
$\Sigma$ is not a submanifold, as is shown by Examples 7.2
and 7.4. Using Theorem 3, we prove that if
$\Aut(M,p_0)$ is discrete at a finitely nondegenerate point
$p_0$, then $\Aut(M,p)$ is also discrete for all $p$ in a
neighborhood of $p_0$ (Theorem 4).

We should mention here that most of the results
described above are optimal, since they fail
to hold  for holomorphically degenerate
hypersurfaces.

Considerable work has been done in the past on
transformation groups of Levi nondegenerate hypersurfaces,
beginning with the seminal paper by Chern--Moser [CM]. (See
also Burns--Shnider [BS] and Webster [W]). Further
contributions were made by the Russian school (see e.g.
the survey papers of Vitushkin [Vi] and Kruzhilin [Kr], as
well as the references therein). Results for higher
codimensional quadratic manifolds were obtained by Tumanov
[T]. We would like to point out that some of our results
are new even for Levi nondegenerate hypersurfaces, since
our approach is quite different, i.e. not based on
the work [CM].

The authors would like to thank Nancy Stanton for several
useful discussions on the subject of this paper.  The second
author would like to thank UCSD for its hospitality during
the preparation of this paper.

\heading \S 1 Jet groups and biholomorphisms of germs of
hypersurfaces\endheading
 Let $k$ be a positive integer and
$
\jet kp=J^k (\bC^N)_p$  the set of
$k$-jets at $p$ of holomorphic mappings from $\bC^N$ to
$\bC^N$ fixing
$p$.  Let
$\gp k= \gp k (\bC^N)$ be the subset of $\jet k 0$
corresponding to holomorphic mappings with nonvanishing
Jacobian determinant at $0$.  $\jet k 0$ can be identified
with the space of holomorphic polynomial  mappings of degree
$\le k$, mapping $0$ to $0$.  The subset $\gp k$ is a
complex Lie group under composition.  The local coordinates
of $\gp k$ can be taken to be the coefficients
$\Lambda=(\L_\a)$ of the polynomials corresponding to the
jets.  The group multiplication then consists of composing
the polynomial mappings and dropping the monomial terms of
degree higher than $k$.

For $p,p' \in \bC^N$, denote by $\scrE_{p,p'}$ the space of
germs of holomorphic mappings $H: (\bC^N,p)
\to (\bC^N,p')$, (i.e. $H(p) = p'$) with Jacobian
determinant of $H$ nonvanishing at $p$. The space
$\scrE_{p,p'}$ is equipped with its natural inductive limit
topology corresponding to uniform convergence on compact
neighborhoods of $p$. That is, a sequence $\{H_j\}
\subset \scrE_{p,p'}$ converges to $H \in
\scrE_{p,p'}$ if there is a compact neighborhood
of $p$ to which all the $H_j$ extend and on which the
$H_j$ converge uniformly to $H$.  If
$H
\in \scrE_{p,p'}$, let $F \in \scrE_{0,0}$ be defined by
$F(Z) = H(Z+p)-p'$.  Then $j_k(F)$, the $k$-jet of $F$ at
$0$, is an element of $\gp k$.  For each
$k$ and $p,p'$ fixed we can define the mapping
$\eta_{p,p'}:
\scrE_{p,p'}\to \gp k$ by
$\eta_{p,p'}(H) = j_k(F)$. In local coordinates $Z$ we have
$\eta_{p,p'}(H)=(\partial^\a_ZH(p))_{1\le |\a| \le k}$. The
mapping
$\eta_{p,p'}$ is continuous;  composition of mappings is
related to group multiplication in $\gp k$ by the identity
$$\eta_{p,p''}(H_2\circ H_1) =
\eta_{p',p''}(H_2)\cdot\eta_{p,p'}(H_1)
\tag 1.1 $$  for any $H_1
\in
\scrE_{p,p'}$ and $H_2 \in
\scrE_{p',p''}$, where $\cdot$ denotes the group
multiplication in $\gp k$. We write $\eta$ for
$\eta_{p,p'}$ when there is no ambiguity.

If $M$ and $M'$ are two real analytic hypersurfaces in
$\bC^N$ and $p$ and $p'$ are points in $M$ and
$M'$ respectively, denote by $\scrF= \scrF(M,p;M',p')$  the
subset of $\scrE_{p,p'}$ consisting of those germs of
mappings which send $M$ into $M'$.  We equip $\scrF$ with
the induced topology.
\proclaim {Theorem 1}  Let $M$ and $M'$ be two real
analytic hypersurfaces in $\bC^N$ which are
$k_0$-nondegenerate at $p$ and $p'$ respectively and let
$\scrF= \scrF(M,p;M',p')$ as above. Then the restriction of
the map
$\eta: \scrE_{p,p'} \to \gp {2k_0}$ to $\scrF$ is
one-to-one; in addition, $\eta(\scrF)$ is a totally real,
closed, real analytic submanifold of $\gp {2k_0}$ (possibly
empty) and $\eta$ is a homeomorphism of $\scrF$ onto
$\eta(\scrF)$. Furthermore,  global defining equations
for the submanifold $\eta(\scrF)$ can be explicitly
constructed from local defining equations for
$M$ and $M'$ near $p$ and $p'$.
\endproclaim

With the notation above, $\scrF(M,p;M,p)$ is the same as
$\Aut(M,p)$, the stability group of $M$ at
$p$.  The algorithm for obtaining defining
equations for $\eta(\scrF)$ described in the proof of
Theorem 1 can be used to calculate $\Aut (M,p)$  and to
determine whether two hypersurfaces are locally equivalent
(see \S7 for examples).

\proclaim {Corollary 1.2} If $M=M'$ and $p=p'$ in  Theorem
1, then $\eta(\scrF)$ is a closed, totally real Lie subgroup
$G(M,p)$ of $\gp {2k_0}$.  Hence the stability group
$\Aut (M,p)$ of $M$ at
$p$ has a natural Lie group structure. In general, for
different $(M,p)$ and $(M',p')$,
$\eta(\scrF)$ is either empty or is a coset of the subgroup
$G(M,p)$.
\endproclaim Combining Theorem 1 and results in [BER2] we
obtain the following.
\proclaim {Corollary 1.3} Let $M$ be a real analytic
connected real hypersurface in
$\bC^N$ which is holomorphically nondegenerate.  Let
$\ell$ be the Levi number of $M$.  Then there is a proper
real analytic subvariety $V \subset M$ such that for any $p
\in M \backslash V$, $\eta$ is a homeomorphism between
$\Aut(M,p)$ and a closed, totally real Lie subgroup of $\gp
{2\ell}$.  Conversely, if $M$ is as above but not
holomorphically nondegenerate, then for any positive
integer $k$ and any $p \in M$, the mapping $\eta: \Aut(M,p)
\to \gp k$ is not injective.
\endproclaim

Note that the complex Lie group $G^k$ can be identified
with a Zariski open subset of the affine space $\bC^\nu$,
for an appropriate choice of $\nu$. A real hypersurface of
$\bC^N$ is called algebraic if it is given by the
vanishing of a real valued polynomial.  For such
hypersurfaces we have the following result.

\proclaim {Theorem 2} If, in addition to the assumptions of
Theorem 1, $M$ and $M'$ are real algebraic hypersurfaces,
then $\eta(\scrF)$ is a real algebraic submanifold of $\gp
{2k_0}$.
\endproclaim

The proofs of Theorems 1 and 2 are given in the next two
sections.
\heading \S 2  First part of the proof of Theorem 1
\endheading For a fixed real hypersurface $M$ and
$p\in M$, we take normal coordinates $Z=(z,w)$ vanishing at
$p$, with $z\in
\bC^n$, $n=N-1$, $w\in \bC$, and assume that $M$ is given by
$w= Q(z,\bar z, \bar w)$, with $Q$ holomorphic near the
origin in $\bC^{2n+1}$ and $Q(z,0,\tau)\equiv
Q(0,\chi,\t)\equiv \t$ (see [CM] or [BJT]). We use analogous
notation for
$(M',p')$; i.e. $Z' = (z',w')$ and so forth. With these
normal coordinates, we may write any $H
\in
\scrF(M,p;M',p')$ in the form $H = (f,g)$, such that the
map is defined by
$z' = f(z,w)$ and $w'=g(z,w)$.  Note that it follows from
the normality of the coordinates that
$g(z,0)\equiv 0$. For each fixed $k$ we choose coordinates
$\Lambda$ in $\gp k$ with $\Lambda = (\lambda_{z^\a
w^j},\mu_{z^\a w^j}), 0< |\a|+j \le k$, such that if $H
=(f,g)
\in \scrF$, then the coordinates of $\eta(H)$ are defined
by $\lambda_{z^\a w^j}=\partial_{z^\a w^j}f(0)$ and
$\mu_{z^\a w^j}=\partial_{z^\a w^j}g(0)$.  We identify an
element in $\gp k$ with its coordinates $\L$.  We shall
denote by $\gp k_0$ the submanifold of $\gp k$ consisting
of those
$\L=(\l,\mu)$ for which $
\mu_{z^\a} = 0$ for all $0 < |\a| \le k$.  It is easily
checked that $\gp k_0$ is actually a subgroup of
$\gp k$ and hence a Lie group. We note also that $\gp k_0$
is stable under complex conjugation, i.e. $\L
\in
\gp k_0$ if and only if $\bar \L \in \gp k_0$.

We associate to $M$ the complex hypersurface $\scrM$ in
$\Bbb C^{2N } $ locally defined near $(p, \bar p)$ by
$$
\scrM=\{(Z,\z): \rho(Z, \z) = 0\},
$$ where $\rho(Z,\bar Z)$ is a real analytic defining
function for
$M$  near
$p$.
Thus, in normal coordinates, $\scrM$ is defined by
$\t = \qb(\chi,z,w)$ for $(z,\chi,w,\t) \in
\bC^{2n+2}$ and similarly for $\scrM'$. If we embed
$\bC^N$ in $\bC^{2N}$ as the totally real plane $\z=\bar
Z$, then $\scrM \cap \bC^N = M$.  Since $H = (f,g)$ maps
$M$ into $M'$,
 it follows from the above that we have for
$(z,w,\chi,\t) \in \scrM$
$$
\gb(\chi,\t) = \qb'(\fb(\chi,\t), f(z,w), g(z,w)),\tag 2.1
$$ where we have used the notation $\bar h(Z) = \overline
{h(\bar Z)}$.

We now introduce the following holomorphic
vector fields which are tangent to
$\scrM$:
$$
\scrL_j = \prt { \chi_j } +\qb_{\chi_j} (\chi,z,w)
\prt {\t}  ,   \ \  j = 1,\ldots n.
\tag 2.2
$$ Note that the
$\scrL_j$ commute with each other.
\proclaim {Lemma 2.3}  Let $M$ and $M'$ be two real
analytic real hypersurfaces in $\bC^N$ which are
$k_0$-nondegenerate at $p$ and $p'$ respectively. Then
there exist
$\bC^N$-valued functions
$\Psi_j(z,
\Lambda)$, $j=0, 1,2,\ldots$, each holomorphic in a
neighborhood of
$0\times
\gp {k_0+j}_0$ in $\bC^n \times \gp {k_0+j}_0$ such that if
$H(z,w) \in \scrF(M,p;M',p')$ with $(\partial^\a
H(0))_{|\a|\le k_0+j} =
\L_0\in
\gp {k_0+j}_0$, then
$$\partial^j_wH(z,0) = \Psi_j(z,\bar \L_0), \ j=0,
1,2,\ldots.$$ Furthermore, we have $\Psi_{0,N}(z,\L)
\equiv 0$, where $\Psi_{0,N}$ is the last component of the
mapping
$\Psi_0$.
\endproclaim
\demo {Proof}Assume $H=(f,g) \in \emmp$.  After applying
the $\scrL_j$ to (2.1) $|\a|$ times, and applying Cramer's
rule after each application, we obtain for
$(z,w,\chi,\t)\in \scrM$ near the origin,
$$ \multline
\qb_{\chi^\a}'(\fb(\chi,\t), f(z,w),g(z,w))=\\
\sum_{|\b|\le |\a|} (\scrL^\b\gb (\chi,\t))
u_{\a,\b}\left((\scrL^\gamma
\fb(\chi,\t))_{|\gamma|\le
|\a|}\right)/\Delta^{2|\a|-1},\endmultline
\tag 2.4$$
where $\Delta =
\Delta(z,w,\chi,\t)=
\det[\scrL_j\fb_k (\chi,\t)]$, and $u_{\a,\b}$
 are universal polynomials, i.e. independent of $M$,
$M'$, and
$H$.  (Note that
$\Delta(0)\not= 0$ since $H$ is assumed to be a
diffeomorphism.)  Since for any germ at $0$ of a
holomorphic function $h(\chi,\t)$ and any multi-index $\b$
we have
$\scrL^\b h(0)= \partial_{\chi^\b}h(0)$, we have
$\scrL^\b g(0) = 0$ by the normality of the coordinates.
Hence the right hand side of (2.4) vanishes at the origin.

   By the assumption that $M'$ is
$k_0$-nondegenerate at $p'$, there exist $n$ multi-indices
$\a^1, \ldots, \a^n$, with $1\le |\a|\le k_0$, such that
$\det[\qb'_{\chi^{\a^j},z_k}(0)]\not= 0$. Hence, by the
implicit function theorem, there exists a unique germ at
$0$ of a holomorphic
 function in $\bC^n\times\bC\times\bC^n$,
$S(\chi',\t',r)$, vanishing at $0$, so that
$X=S(\chi',\t',r)$ solves the system of equations
$$\qb_{\chi^{\a^j}}'(\chi',X,Q'(X,\chi',\t'))=r_j,
\ \ j=1,\ldots,n,  \tag 2.5 $$
with $r=(r_1,\ldots,r_n)$. Note that $S$ depends only on
$M', p'$. We take in (2.4) $\chi = 0$ and $\t=w$, so that
$z,w$ are free, since $(z,w,0,w) \in \scrM$.  Denote by
$R^H_\a(z,w)$ the right hand side of (2.4) after making
this substitution. Hence we have
$$f(z,w) =
S\left(\fb(0,w),\gb(0,w),(R^H_{\a^j}(z,w))_{1\le j\le
n}\right).\tag 2.6$$

It follows from the explicit definition of
$R^H_\a(z,w)$ that for any $\a, |\a|\le k_0$, and
$j=0,1,2,\ldots$ there is a function
$J_{\a,j}(z,\L)$, holomorphic in a neighborhood of $0\times
\gp {k_0+j}_0$ in $\bC^n\times\gp{k_0+j}_0$, such that
$$\partial^j_wR^H_\a(z,0)=J_{\a,j} (z,(\partial^\gamma \bar
H(0))_{|\gamma|\le k_0+j}), \quad j=0,1,\ldots.
\tag 2.7$$ The functions $J_{\a,j}(z,\L)$ depend only on
$M,p$ (but not on $M'$ or $H$).  The conclusion of the
lemma for the $f$ components now follows by differentiating
(2.6) with respect to
$w$ $j$ times and then making use of (2.7).  The conclusion
for
$g$ then follows by differentiating the complex conjugate of
(2.1) in $w$ and substituting the value of
$\partial^j_w f(z,0)$ found above. This completes the proof
of Lemma 2.3. $\square$
\enddemo
\proclaim {Lemma 2.8} Let $M, p, M', p'$ be as in Lemma
2.3.  Then there exists a $\bC^N$-valued function
$\Phi(z,\chi,\L)$,  holomorphic in a neighborhood of
$0\times 0\times\gp {2k_0}_0$ in $\bC^n \times \bC^n
\times \gp {2k_0}_0$, such that for $H \in \emmp$ with
$(\partial^\a H(0))_{|\a|\le 2k_0} = \L_0$, we have
$$H(z,Q(z,\chi,0)) \equiv \Phi(z,\chi,\L_0).$$
\endproclaim
\demo {Proof} We begin with (2.4) in which we take
$\t=0$ and $w=Q(z,\chi,0)$, making $z$ and $\chi$ free
variables.  Denote by $\scrR^H_\a(z,\chi)$ the right hand
side of (2.4) with this substitution. Note that
$\scrR^H_\a(0,0)=0$.  Let $S(\chi',\t',r)$ be the solution
of (2.5) as in the proof of Lemma 2.3. By the implicit
function theorem, as in the proof of Lemma 2.3, we have
$$f(z,Q(z,\chi,0)) =
S\left(\fb(\chi,0),0,(\scrR^H_{\a^j}(z,\chi))_{1\le
j\le n}\right).\tag 2.9$$   By
taking complex conjugates in the conclusion of Lemma 2.3
for $0\le j\le k_0$ and substituting in the expression for
$\scrR^H_\a(z,\chi)$ it is easy to see that for any
$|\a|\le k_0$ there exist  functions
$\scrJ_\a(z,\chi,\L)$ holomorphic in a neighborhood of
$0\times0\times \gpo$ in
$\bC^n\times\bC^n\times\gpo$ satisfying
$$\scrR^H_\a(z,\chi)=\scrJ_\a(z,\chi, (\partial^\b
H(0))_{|\b|\le 2k_0}). \tag 2.10$$ It is important to note
here that the functions
$\scrJ_\a$ depend only on $M,p$, $M',p'$ and not on
$H$. The conclusion of Lemma 2.8 for the
$f$ components of
$H$ then follows by substituting (2.10) in (2.9) and by
using the expression for
$\fb(\chi,0)$ given by Lemma 2.3.  As in the proof of Lemma
2.3, the expression for the $g$ component follows from the
conclusion for the $f$ components by using (2.1).
  $\square$
\enddemo

\proclaim {Proposition 2.11}Let $M$ and $M'$ be two real
analytic real hypersurfaces in $\bC^N$ which are
$k_0$-nondegenerate at $p$ and $p'$ respectively. There
exists a $\bC^N$-valued function
$F(z,t,\L)$ holomorphic in a neighborhood of
$0\times 0\times \gp {2k_0}_0 $ in
$\bC^n\times\bC\times\gp {2k_0}_0$ and a germ at $0$ of a
nontrivial holomorpic function $B(z)$,  such that for a
fixed
$\L_0
\in
\gp {2k_0}_0$ there exists $H \in \emmp$ with
$$(\partial^\a H(0))_{|\a|\le 2k_0}=\L_0 \tag 2.12$$
 if and only if all of the following hold:
\roster
\item "(i)" $(z,w) \mapsto F(z,{w\over B(z)},\L_0)$ extends
to a function $K_{\L_0}(z,w)$ holomorphic in a full
neighborhood of $0$ in $\bC^N$.
\item "(ii)" $(\partial^\a K_{\L_0}(0))_{|\a|\le 2k_0}=\L_0
$.
\item "(iii)" $K_{\L_0}(M) \subset M'$.
\endroster If {\rm (i), (ii), (iii)} hold, then the unique
mapping in
$\emmp$ satisfying (2.12) is given by
$H(Z)=K_{\L_0}(Z)$.
\endproclaim
\demo {Proof} We start with Lemma 2.8. From the
$k_0$-nondegeneracy, we necessarily have
$Q_{\chi_1}(z,0,0)\not\equiv 0$ and we set
$$A(z) = Q_{\chi_1}(z,0,0).
\tag 2.13 $$ We write $\chi=(\chi_1,\chi')$; we shall solve
the equation
$$ w=Q(z,(\chi_1,0),0)\tag 2.14$$ for $\chi_1$ as a
function of $(z,w)$ and analyze the solution as $z$ and $w$
approach $0$. We have
$$Q(z,(\co,0),0)= \sum_{j=1}^{\infty}A_j(z)\co^j,
\tag 2.15$$ with $A_1(z) = A(z)$ and $A_j(0) = 0$, $j =
1,\ldots$. Dividing (2.14) by $A(z)^2$, we obtain
$$ {w\over A(z)^2} = {\co\over A(z)} +
\sum_{j=2}^\infty A_j(z){\co^j\over [A(z)]^2}. $$ We set
$C_j(z) = A_j(z)A(z)^{j-2}$, $j \ge 2$, and let
$$\psi(z,t)= t + \sum_{j=2}^\infty v_j(z)t^j \tag 2.16$$
 be the solution in $u$ given by the implicit function
theorem of the equation
$t = u + \sum_{j=2}^\infty C_j(z)u^j,$ with $\psi(0,0)=0$.
The functions $\psi$ and $v_j$ are then holomorphic at $0$
and $v_j(0) = 0$.  A solution for $\co$ in (2.14) is then
given by
$$\co = \theta(z,w) = A(z) \psi\left(z, {w \over
A(z)^2}\right).\tag 2.17$$ The function $\theta(z,w)$ is
holomorphic in an open set in $\bC^{n+1}$ having the origin
as a limit point.

Now define $F$ by
$$F(z,t,\L) = \Phi(z,(A(z)\psi(z,t),0),\L),\tag2.18$$ where
$\Phi$ is given by Lemma 2.8, and let $B(z) = A(z)^2$,
with $A(z)$ given by (2.13). Then (i) follows from Lemma
2.8.  The rest of the proof of the proposition is now easy
and is left to the reader.
$\square$
\enddemo
\heading \S 3 End of the proof of Theorem 1; proof of
Theorem 2\endheading In this section we shall give
equivalent conditions for (i),(ii), (iii) of  Proposition
2.11, which will imply Theorem 1.
\proclaim {Proposition 3.1} Under the hypotheses and
notation of Proposition 2.11, there exists a function
$K(Z, \L)$, holomorphic in a neighborhood of
$0\times \gpo$ in $\bC^N
\times\gpo$ such that (i) holds for a fixed
$\L_0 \in \gpo$ if and only if
$F(z, {w\over B(z)},\L_0) \equiv K(z,w,\L_0)$. Furthermore,
we have the following equivalences.
\roster
\item "(a)" There exist functions $c_j$, $j=1,2,\ldots$,
holomorphic in $\gp {2k_0}_0$ such that (i) holds if and
only if  $c_j(\L_0)=0$, $j=1,2,\ldots$.
\item "(b)" There exist functions
$d_j$,
$1\le j\le J$,   holomorphic in $\gp {2k_0}_0$ such that
if (i) is satisfied then (ii) holds if and only if
$d_j(\L_0)=0$,
$1\le j\le J$.
\item "(c)" There exist functions
$e_j$,
$j=1,2,\ldots$, holomorphic in $\gpo\times\gpo$ such that
if (i) is satisfied then (iii) holds if and only if
$e_j(\L_0,
\bar
\L_0)=0$,
$j=1,2,\ldots$.
\endroster
\endproclaim
\demo {Proof} Recall that $F(z,t,\L)$ is holomorphic in a
neighborhood of $0\times 0\times\gpo$ in
$\bC^n\times\bC\times\gpo$. Hence we can write
$$F(z,t,\L) = \sum_{\a,j} F_{\a j}(\L) z^\a t^j,
\tag 3.2$$ with $F_{\a,j}$ holomorphic in $\gpo$. For each
compact subset  $L \subset \gpo$  there exists $C>0$
such that the series $(z,t)\mapsto \sum_{\a,j} F_{\a j}(\L)
z^\a t^j$ converges uniformly for $|z|, |t|\le C$ and for
each fixed $\L \in L$. For $|z|\le C$ and
$|{w\over B(z)}| \le C$ we have
$$ F\left(z,\wob,\L\right)= \sum_{j=0}^\infty {F_j(z,\L)
\over B(z)^j}w^j \tag 3.3$$ with $F_j(z,\L)=\sum_\a
F_{\a,j}(\L)z^\a$. After a linear change of holomorphic
coordinates if necessary, and putting $z=(z_1,z')$, we may
assume, by using the Weierstrass Preparation Theorem, that
$$B(z)^j = U_j(z)[z_1^{K_j} +
\sum_{p=0}^{K_j-1}a_{jp}(z')z_1^p],$$ with $U_j(0)\not=0$
and $a_{jp}(0)=0$. By the Weierstrass Division Theorem we
have the unique decomposition
$$F_j(z,\L)=Q_j(z,\L)B(z)^j+
\sum_{p=0}^{K_j-1}r_{jp}(z',\L)z_1^p, \tag 3.4
$$ where $Q_j(z,\L)$ and $r_{jp}(z',\L)$ are holomorphic in
a neighborhood of $0\times\gpo$ in
$\bC^n\times\gpo$. Moreover, there is a constant $C_1>0$
(see e.g. [H]) such that
$$
\sup_{|z|\leq\delta}|Q_j(z,\Lambda)|\leq C_1^j
\sup_{|z|\leq\delta}|F_j(z,\Lambda)|,
$$ where $\delta$ is some number, $0<\delta<C$. This,
together with (3.3) implies that for a given
$\L_0\in \gpo$,  (i) holds if and only if
$z'\mapsto r_{jp}(z',\L_0)$ vanishes identically for all
$j,p$.

The first statement in Proposition 3.1 follows by taking
$$K(z,w,\L) = \sum_j Q_j(z,\L)w^j.\tag 3.5$$
 Taking the Taylor
expansion $r_{jp}(z',\L)= \sum _\a c_{jp\a}(\L){z'}^\a$, we
see that (i) is equivalent to the vanishing of all the
$c_{jp\a}(\L_0)$, which yields (a).  Since $K_{\L_0}(Z)
\equiv K(Z,\L_0)$, (b) is proved by taking $d_j(\L)$ as the
components of
$(\partial_Z^\a K(0,\L))_{|\a|\le 2k_0}-\L$. For (c), we
note that (iii) is equivalent to
$$\rho'(K(z,w,\L_0), \bar K(\chi,
\qb(\chi,z,w),\bar \L_0)) \equiv 0, \tag 3.6$$ where
$\rho'$ is a defining function for $M'$.  We obtain the
conclusion in (c) by expanding the left hand side of (3.6)
as a series in
$z,w,\chi$ with coefficients which are holomorphic functions
of
$\L_0,\bar \L_0$.  $\square$
\enddemo
\demo {Proof of Theorem 1} It follows immediately from
Proposition 2.11 that the mapping
$\eta:\scrE_{p,p'}\to \gpo$ is one-to-one.  From
Proposition 3.1 we conclude that $\eta(\scrF)$ is a closed
subset of $\gpo$, since it is defined by the vanishing of a
set of real analytic functions. The continuity of $\eta$ is
clear. To see that $\eta$ is a homeomorphism onto its
image, we observe that if
$\L_j
\in
\eta(\scrF)$ converges to $\L_0$, then $\L_0 \in
\eta(\scrF)$.  Also, by Proposition 3.1,
$\eta^{-1}(\L_j)=K(Z,\L_j)$ and
$\eta^{-1}(\L_0)=K(Z,\L_0)$. Since $K(Z,\L)$ is holomorphic
in a neighborhood of $(0,\L_0)$, it is obvious that
$\eta^{-1}(\L_j))$ converges to
$\eta^{-1}(\L_0)$.  This proves that $\eta$ is a
homeomorphism onto $\eta(\scrF)$.

To prove that $\eta(\scrF)$ is actually a manifold, we
consider first the case where $M=M'$ and $p=p'$, so that
$\emmp=\Aut(M,p)$. In this case,
$\eta(\scrF)$ is actually a subgroup of $\gpo$ and
$\eta$ is a group homomorphism.  Since $\eta(\scrF)$ is
closed, it is a Lie subgroup (see e.g. [Va]), and hence a
real analytic manifold in the induced topology.  In the
general case it is easy to check that
$\eta(\emmp)$ is a coset of $\eta(\Aut(M,p))$ in the group
$\gpo$ and therefore is also a manifold, which is either
empty or homeomorphic to $\Aut(M,p)$.

Finally, the proof that the parameter set
$\eta(\emmp)$ is a totally real manifold follows easily from
the following lemma.

\proclaim {Lemma 3.7} There exists a holomorphic mapping
$T:\gpo\to \gpo$ such that for every $\L \in
\eta(\emmp)$,
$$ \L =T(\bar \L). $$
\endproclaim
\demo {Proof of Lemma 3.7}  We begin with Lemma 2.3, from
which we conclude that for any $j =
0,1,\ldots, $, any $n$-multi-index
$\b$, and any $H \in \scrF(M,p;M',p')$ we have
$$\partial_z^\b\partial_w^j H(0) = \partial_z^\b
\Psi_j\left(0,\left(\partial^\a \bar
H(0)\right)_{|\a|\le k_0+j }\right),\tag 3.8
$$ where the $\Psi_j(z,\L)$ are given by Lemma 2.3.
Let $K(Z,\l)$ be the mapping given by Proposition 3.1.
 For any $(n+1)$-multi-index $\a$, let $\scrK_\a$ be the
holomorphic function on $\gpo$ given by $\scrK_\a(\L)
= \partial_Z^\a K(0,\L)$.  It follows from Propositions
2.11 and 3.1 that the following holds for any $H \in
\scrF(M,p;M',p')$ and any $(n+1)$-multi-index $\a$.
$$ \partial^\a H(0) = \scrK_\a
\left((\partial^\b H(0))_{|\b| \le 2k_0}\right). \tag 3.9$$
Lemma 3.7 follows by making use of (3.8) and the conjugate
of (3.9) for $|\a|\le 3k_0$.
$\square$
\enddemo

The proof of Theorem 1 is now complete.
\qed
\enddemo

\demo {Proof of Theorem 2}
   It suffices to go through the
arguments used in the proof of Theorem 1 and to observe
that  all functions and mappings that appear are algebraic,
since they are obtained from polynomials by
differentiation, the implicit function theorem, and the
Weierstrass Preparation Theorem (cf. [BR] and [BER1]  for
similar arguments). Moreover, since the defining functions
of $\eta(\scrF)$ given by Proposition 3.1 are global,
standard algebraic geometry implies that they are in fact
rational.  Hence $\eta(\scrF)$ is a real algebraic
submanifold in $\gpo$. We leave the details to the reader.
\qed
\enddemo
\heading \S 4 Jet bundles and dependence on base points
\endheading In this section we shall generalize the
results of \S 1 to the case where $p$ and $p'$ are varying
points in $M$ and
$M'$ respectively.  We begin with some notation and
definitions.

If $X$ and $Y$ are two complex manifolds and $k$ a positive
integer, we denote by  $\jkxy$ the complex manifold of
$k$-jets of germs of holomorphic mappings from $X$ to $Y$,
i.e.
$$ \jkxy = \bigcup_{x\in X,y\in Y}\jkxy_{(x,y)} $$ where $
\jkxy_{(x,y)}$ denotes the $k$-jets of germs at $x$ of
holomorphic mappings from $X$ to $Y$ and taking $x$ to
$y$. (See e.g. [M], [GG].) With this notation,
$ J^k(X,X)_{(x,x)}$ is the same as $J^k(X)_x$ introduced in
\S 1 with $X=\bC^N$.

Denote by $\exy$ the set of germs of holomorphic mappings
from
$X$ to $Y$, and similarly by $\exy_{(x,y)}$ those germs
at
$x$ mapping
$x$ to
$y$.  We equip $\exy_{(x,y)}$ with the natural
inductive limit topology used in the previous sections.
(We do not consider a topology on all of $E(X,Y)$.)

For every $k$ there is a canonical mapping
$\sigma_k:\exy\to\jkxy$.  Note that
$\sigma_k|_{\exy_{(p,p')}}$ is the same as the mapping
$\eta_{p,p'}$ defined in \S1 with $X=Y=\bC^N$ (and
hence is continuous.) If $\dim_\bC X=
\dim_\bC Y$ then we denote by
$\gkxy$ the open complex submanifold of $\jkxy$ given by
those jets which are locally invertible. Similarly, we
denote by
$\sexy$ the  subset of $\exy$ consisting of the
invertible germs.  It is clear that the restriction of
$\sigma_k$ maps
$\sexy$ to $\gkxy$.

If $M\subset X$ and $M'\subset Y$ are real analytic
submanifolds, we let $\emxy$ be the set of germs
$H_p\in\exy$ with
$p\in M$
  which map a neighborhood of $p$ in $M$ into $M'$.
Similarly, we denote by $\semxy$ those germs in $\emxy$
which are invertible.  (Note that with this notation
$$\semxy_{(p,p')}=\scrF(M,p;M',p')$$  in the
notation used in the previous sections with
$X=Y=\bC^N$.)
\proclaim {Theorem 3} Let $M\subset X$ and
$M'\subset Y$ be two real analytic hypersurfaces with
$\dim X =\dim Y$, and
$k_0$ a positive integer.  Suppose that $M$ and $M'$ are
both at most
$k_0$-nondegenerate at every point. Then the mapping
$$\sg_{2k_0}: \semxy\to \gp {2k_0}(X,Y)$$ is one-to-one
onto its image $\Sigma$.  Furthermore,
$\Sigma$ is a real analytic subset of $\gp {2k_0} (X,Y)$,
possibly empty, and each fiber $\Sigma\cap\gp {2k_0}
(X,Y)_{(p,p')}$, with $p\in M$, $p'\in M'$, is a real
analytic submanifold homeomorphic to
$\scrF(M,p;M',p')$.
\endproclaim

\remark {Remarks}
\roster
\item"(i)" If $M\subset X$ and $M'\subset Y$ are real
analytic submanifolds, we write
$$\gp k (X,Y)_{(M,M')} =\bigcup_{(p,p')\in M\times M'} \gp
k (X,Y)_{(p,p')}.
$$ It is easy to check that $\gp k (X,Y)_{(M,M')}$ is a
real analytic submanifold of $\gp k (X,Y)$.  We observe
here that in Theorem 3 above, the real analytic set
$\Sigma$ is actually a subset of $\gp {2k_0}
(X,Y)_{(M,M')}$.
\item"(ii)" Under the assumptions of Theorem 3,   if
$M$ and
$M'$ are
$k$-nondegenerate at $p\in M$ and $p'\in M'$, respectively,
for some $k<k_0$, then $\semxy_{(p,p')}$ can be embedded
into the smaller space $G^{2k}(X,Y)_{(p,p')}$ according to
Theorem 1. In the special case where $M$ and
$M'$ are both $k_0$-nondegenerate at every point (hence
$\ell(M)=\ell(M')=k_0$), Theorem 3 is a
 generalization of Theorem 1 with $p$ and $p'$ allowed
to vary.
\item"(iii)" It follows from Theorem 3 that the
set $\semxy$ may be given the topology obtained by
identifying it with
$\Sigma$; this topology is compatible with the
topology of each fiber
$\semxy_{(p,p')}$.
\endroster
\endremark

We write
$$\gp k (X)= \bigcup_{p\in X}
 \gp k (X,X)_{(p,p)}$$ for the fiber bundle of invertible
jets fixing a point. Note that each fiber is isomorphic to
the Lie group $\gp k (\bC^N)$, with $N = \dim X$.   If
$M\subset X $ is a real analytic hypersurface, we let
$\scrA(M)=
\cup_{p\in M}
\Aut(M,p)$.
  We have the following.
\proclaim {Corollary 4.1} Let $M\subset X$ be a real
analytic hypersurface and $k_0$ a positive integer.  Assume
that $M$ is at most
$k_0$-nondegenerate at all points.  Then the mapping
$$\sg_{2k_0}: \scrA(M)\to \gp {2k_0}(X)$$ is
one-to-one onto its image $\Sigma$.  Furthermore,
$\Sigma$ is a real analytic subset of $\gp {2k_0} (X)$,
 and each fiber $\Sigma\cap\gp {2k_0} (X)_{p}$, with $p\in
M$, is a closed real Lie subgroup of the complex Lie group
$\gp {2k_0} (X)_{p}$ homeomorphic to $\Aut(M,p)$.
\endproclaim As in Remark (i) above, the real
analytic set $\Sigma$ of  Corollary 4.1 is in fact a subset
of the subbundle
$\gp {2k_0} (X)_M$.
\demo {Proof of Theorem 3}  Since the conclusion of the
theorem is local, we pick $\po \in M$, $\po' \in M'$ and
assume $X=Y=\bC^N$. The following lemma, whose proof easily
follows from the construction of normal coordinates at a
given point (see [CM] [BJT]), shows that normal coordinates
can be chosen to depend real analytically on the central
point of a real analytic hypersurface. Let
$\scrM$ be the complexification of $M$ as in \S 2.

\proclaim {Lemma 4.2}  Let $M$ be a real analytic
hypersurface in $\bC^N$ and $\po \in M$. Then there is a
mapping
$(z(Z, p,\bar p),w(Z, p,\bar p))$ from a neighborhood of
$(\po,(\po,\bar\po))$ in
$\bC^N\times\scrM$  into $\bC^n\times\bC$ such that
$(z(p, p,\bar p),w(p, p,\bar p))=0$ for all $p\in M$ near
$\po$, $Z\mapsto (z(Z, p,\bar p),w(Z, p,\bar
p))$ is a change of coordinates near $p$, and there
exists a function $Q(z,\chi,\t;p,q)$ holomorphic in a
neighborhood of
$(0,0,0,(\po,\bar\po))$ in $\bC^{2n+1}\times \scrM$ such that
$$Q(z,0,\t;p,q)\equiv Q(0,\chi,\t;p,q)\equiv \t$$ and such
that
$$ w-Q(z,\bar z,\bar w;p,\bar p)=0\tag4.3
$$  defines $M$ in a neighborhood of $p$ in $\bC^N$.
\endproclaim

We return to the proof of Theorem 3. For any $p\in M$ near
$p_0$ and
$p'\in M'$ near
$\po'$, we choose normal coordinates as in Lemma 4.2 and
write the defining equations of $M$ and $M'$ as in (4.3).
The conclusion of Theorem 1 still holds and
its proof is the same if we replace $G^{2k_0}_0$ by
$G^{2k}_0$ for any
$k>k_0$ (taking the corresponding mapping $\eta$).  By
following the real analytic dependence on the parameters
$p$ and $p'$ in (4.3) and the corresponding defining
equation for $M'$, a detailed inspection of  the proof of
Theorem 1 shows that the defining functions
$c_j$,
$d_j$ and
$e_j$ for the fiber at $(p,p')$ of the image $\Sigma$,
given by Proposition 3.1, vary real analytically with $p$
and $p'$.  From this, the Theorem easily follows.
$\square$
\enddemo

\heading \S 5 Discreteness of $\Aut (M,p)$\endheading

In this section we shall show that if $\Aut (M,\po)$ is
discrete at a finitely nondegenerate point $\po$, the same
is true for $p \in M$ near $\po$.  The main tool in proving
this result is Corollary 4.1 which describes the dependence
of $\Aut (M,p)$ on the point $p$.

\proclaim {Theorem 4} Let $M$ be a real analytic
hypersurface in $\bC^N$ finitely nondegenerate at
$\po$. If $\Aut(M,\po)$ is a discrete group, then
$\Aut(M,p)$ is also discrete for all $p$ in a neighborhood
of $\po$ in $M$.  Equivalently, if $
\hol_0(M,\po)=\{0\}$ then $
\hol_0(M,p)=\{0\}$ for all $p$ in a neighborhood of $\po$
in $M$.
\endproclaim

By using Corollary 4.1 we may reduce the proof of Theorem 4
to that of Lemma 5.1 below. The last statement of
the theorem follows from the first part, since $\hol_0(M,p)$
is the Lie algebra of $\Aut(M,p)$.

The following could be proved using the ``no small
subgroups" property of Lie groups (see e.g. [MZ]), but
we shall give a self-contained proof for the
convenience of the reader.

\proclaim {Lemma 5.1} Let $G$ be a real Lie group and
$c_j(g,x)$, $j = 1, 2, \ldots$,  continuous functions in
$G\times U$, where $U$ is an open subset of $\bR^q$.
Assume that for every $x \in U$ the set $\Sigma_x:=\{g\in
G: c_j(g,x)=0
\}$ is a closed subgroup of $G$. If
$\Sigma_{x_0}$ is discrete for some $\xo\in U$, then
$\Sigma_{x}$ is also discrete for all $x$ in a neighborhood
of $x$ in $U$.
\endproclaim

\demo {Proof of Lemma 5.1} The proof is by contradiction.
Assume that there exists a sequence $x_j\in U$ converging to
$\xo$ such that for every $j$, $ \Sigma_{x_j}$ is not
discrete i.e.
$\dim\Sigma_{x_j} \ge 1$. We shall show that this implies
that the identity, ${\bold 1}$, is not an isolated point in
$\Sigma_{\xo}$, which would contradict the assumption of
the lemma.

 Let $V$ be a sufficiently small neighborhood of the
identity in
$G$ such that there is a diffeomorphic local isomorphism of
$V$ with a neighborhood of the identity in a closed
subgroup of $GL(k, \bC)$, for some $k$, which is possible
by Ado's Theorem (see e.g. [Va]). We  pull back to $V$ the
Euclidean metric on $GL(k,\bC)$ and denote by
$B_r({\bold 1})$ the ball in $V$ of radius $r$ around
${\bold 1}$ for $r$ sufficiently small. We denote by
$\Sigma_{x}^0$ the identity component of $\Sigma_x$. We
claim that there exists
$\delta>0$ such that
$$
\Sigma_{x_j}^0\cap \partial B_\delta({\bold
1})\neq\emptyset.\tag5.2
$$

Since $V$ is diffeomorphically isomorphic to a neighborhood
of the identity in a closed subgroup of
$GL(k,\bC)$, the claim is a consequence of the following
lemma.

\proclaim {Lemma 5.3} Let $k$ be a positive integer.
There exists
$\epsilon>0$ such that for all real subgroups $G'$ of
$GL(k,\bC)$ of positive dimension,
$$ G'\cap\{m\in GL(k,\bC)\: ||m-\bold
1||=\epsilon\}\neq\emptyset,
$$ where $||\cdot||$ is the Euclidean norm in the space of
$k\times k$ matrices with complex coefficients.
 \endproclaim

\demo {Proof of Lemma 5.3} It suffices to assume that
$G'$ is a connected one parameter group $\{\exp tA\: t\in
\bR\}$ with $A$ a nonzero $k\times k$ matrix.  If
$G'$ is unbounded the conclusion will clearly hold for any
$\epsilon > 0$.  Hence we are reduced to the case where
$G'$ is compact, i.e. the matrix $A$ is diagonalizable with
all eigenvalues purely imaginary. We assume this and write
$$A = B\left(\matrix i\b_1&0&\hdots&0\\0&i\b_2&\hdots&0\\
\vdots&\vdots&\ddots &\vdots\\
0&0&\hdots&i\b_k\endmatrix\right)B^{-1}
\tag 5.4$$ with $\b_1,\ldots,\b_k$ real and not all 0 and $B
\in GL(k,\bC)$. We have
$$ \exp tA = B(\exp tD) B^{-1},$$ where $D$ is the diagonal
matrix in (5.4). We note that for any
$k\times k$ matrix
$C$,
$$ |\tr\  C| \le k^{1/2} ||C||. \tag 5.5
$$ Since $$\tr\ (\exp tA-{\bold 1}) = \sum_{j=1}^k
e^{it\b_j}-k, $$ it follows from (5.5) that
$$ |\sum_{j=1}^k e^{it\b_j}-k| \le k^{1/2} ||\exp tA-{\bold
1}||.
\tag 5.6$$ Without loss of generality, we may assume
$\b_1\not= 0$ and choose $t_0$ such that
$e^{it_0\b_1}=-1$. We conclude from (5.6) that
$||\exp t_0A-{\bold 1}||\ge 2/k^{1/2}$. The conclusion of
the lemma follows from the connectedness of $G'$ with
$\epsilon = 2/k^{1/2}$.\qed
\enddemo

We return now to the proof of Lemma 5.1. By connectedness
of $\Sigma_{x_j}^0$ and Claim (5.2), we deduce that for
every $\delta'$, $0 <
\delta' <\delta$, and every $j$ there exists $g_j
\in \Sigma_{x_j}^0\cap \partial B_{\delta'}(\bold 1)$.  By
going to a subsequence if necessary, we may assume that the
sequence $g_j$ converges to an element $g_{\delta'} \in
\partial B_{\delta'}(\bold 1)$.  Since the functions
$c_l(g,x)$ are continuous, we conclude that
$c_l(g_{\delta'},\xo)=0$ for all $l$.  Hence
$g_{\delta'}\in \Sigma_\xo$.  Since $\delta'$ as above is
arbitrary, this shows that $\bold 1$ is not an isolated
point of
$\Sigma_\xo$.
  \qed
\enddemo

\heading\S 6 An algorithm to compute mappings between
hypersurfaces; formal mappings
\endheading

The proof of Theorem 1  actually gives an
algorithm to construct  the defining equations of the
manifold
$$
\Sigma_{p,p'}=\eta(\scrF(M,p;,M',p'))
$$  from  defining equations of $M$ and $M'$ near $p$ and
$p'$. Moreover, for each $\Lambda\in G^{2k_0}_0$ the
algorithm constructs a mapping which is the unique
biholomorphic mapping $H$ sending $(M,p)$ into
$(M',p')$ with
$\eta(H)=\Lambda$ provided $\Lambda\in\Sigma_{p,p'}$. We
summarize here this algorithm, using the notation of \S2 and
\S3.

{\bf Step 1.} From the construction given in the proof of
Lemma 2.3 we obtain the functions $\Psi_j(z,\bar\Lambda)$
for $j=0,1,...,k_0$. This determines the candidate for
$\partial^j_w H(z,0)$.

{\bf Step 2.} From the construction given in the proof of
Lemma 2.8 we obtain the function $\Phi(z,\chi,\Lambda)$
which determines $H(z,Q(z,\chi,0))$.

{\bf Step 3.} Let $A(z)$ be the function defined by (2.13),
and $\chi_1=\theta(z,w)$ the solution of (2.14) given
by (2.17).  From the construction given in the proof of
Proposition 2.11 we obtain the holomorphic function
$F(z,t,\Lambda)$ defined by (2.18), with the property that
$$ H(z,w)=F\left(z,\frac{w}{A(z)^2},\Lambda\right)
$$ provided that $H(z,w)$ is a biholomorphic mapping taking
$(M,p)$ into $(M',p')$ with $\Lambda=\eta(H)$.

{\bf Step 4.} From the construction given in the proof of
Proposition 3.1 we obtain the holomorphic function
$K(z,w,\Lambda)$ given by (3.5), which is the candidate for
$H(z,w)$ if the mapping exists.

{\bf Step 5.} To construct the defining equations for
$\Sigma_{p,p'}$, we expand $F(z,{w \over A(z)^2},\Lambda)$
as a series in $w$ and follow the proof of Proposition 3.1
to find the functions $c_j$, $d_j$, and $e_j$.

All the steps above are constructive since the
implicit function theorem and the Weierstrass preparation
theorem can be implemented by successive iterations.

As a consequence of this explicit construction, we  obtain
the following result on the convergence of formal mappings.
Let $H$ be a formal holomorphic mapping
$(\bC^N,0)\to (\bC^N,0)$, i.e. the components of $H$ are
formal series in $Z$ with vanishing constant term. We
assume that $H$ is invertible, i.e. the formal Jacobian of
$H$ at 0 is invertible. Let $M$ and $M'$ be real analytic
hypersurfaces in $\bC^N$; for simplicity, we assume that
$0\in M$ and $0\in M'$. We say that the formal mapping $H$
maps
$(M,0)$ into $(M',0)$, and write $H(M) \subset M'$, if
$$
\rho'(H(Z),\bar H(\zeta))=a(Z,\zeta)\rho(Z,\zeta)
$$ as formal power series, where $\rho(Z,\bar Z)=0$,
$\rho'(Z,\bar Z)=0$ are defining equations for
$(M,0)$ and $(M',0)$ respectively, and $a(Z,\z)$ is a
formal power series in $(Z,\z)$.

\proclaim {Theorem 5} Let $M$ and $M'$ be real analytic
hypersurfaces in $\bC^N$, both containing $0$, and assume
that $M$ is finitely nondegenerate at $0$. Let
$H
\:(\bC^N,0)\to (\bC^N,0)$ be a formal invertible
holomorphic mapping with $H(M)\subset M'$. Then $H$ is a
germ of a biholomorphism at 0 mapping $M$ into $M'$. In
particular, $(M,0)$ and $(M',0)$  are formally equivalent
if and only if they are biholomorphically equivalent.
\endproclaim

\remark{Remark} The following can be considered a converse
of Theorem 5. If
$M$ is holomorphically degenerate, with $0\in M$, then
there is a formal invertible holomorphic map
$H\:(\bC^N,0)\to(\bC^N,0)$, which is not convergent, such
that $H(M)\subset M'$. To see this, take a germ of a
holomorphic vector field $X$ at 0 tangent to $M$ and
multiply it by a formal, non-convergent power series
$D(Z)$ with no constant term. The formal flow of the formal
vector field $D(Z)X$ gives a $1$-parameter group of formal,
non-convergent, invertible holomorphic mappings preserving
$(M,0)$. We leave the details to the reader. We should
mention here that we do not know of any example of two
germs of real analytic hypersurfaces which are formally
equivalent but not biholomorphically equivalent. There are,
however, examples of $N$ dimensional real analytic
submanifolds in $\bC^N$ which are formally equivalent but
not biholomorphically equivalent. (For this see Gong [G].)
\endremark

\demo{Proof of Theorem 5} Let $k_0$ be the integer such that
$M$ is
$k_0$-nondegenerate at
$0$. It follows easily from the formal equivalence of $M$
and $M'$ at 0 that $M'$ is also
$k_0$-nondegenerate at 0.  If $(z,w)$ are normal
coordinates for $M$ at $0$, we write
$$H(z,w)=\sum_{j=0}^\infty H^j(z)w^j,\tag6.1
$$ where $H^j(z)\in(\bC[[z]])^N$ (the $N$-tuples of formal
power series in
$z$). From the construction in Step 1 above, we have
$$ H^j(z)=\frac{1}{j!}\Psi_j(z,\Lambda_0),\quad j=0,1,...,
$$ where $(\partial ^\alpha H(0))_{|\alpha|\leq
k_0+j}=\Lambda_0\in G^{k_0+j}_0$. This implies that all the
formal power series $H^j(z)$ are in fact convergent.
Similarly, from Step 2 we conclude that the formal power
series
$(z,\chi)\mapsto H(z,Q(z,\chi,0))$ is convergent. It
remains to show that $H(z,w)$ is convergent.

We denote by $\Cal N$ the ring of germs at $0$ of
meromorphic functions $m(z)$ for which $A(z)^km(z)$ is
holomorphic for some integer $k$. We consider $H(z,w)$,
given by (6.1), as an element of $(\Cal N[[w]])^N$, the
$N$-tuples of formal series in $w$ with coefficients in
$\Cal N$. In
$H(z,Q(z,\chi,0))$ we substitute $\chi=(\chi_1,0,...,0)$
with $\chi_1=\theta(z,w)$ as in Step 3 and recover
$H(z,w)$. By Step 3, $H(z,w)$ must coincide, as an element
of
$(\Cal N[[w]])^N$, with $F(z,w/A(z)^2,\Lambda_0)$. By
identifying the coefficients of $w^j$, we conclude that the
coefficients of $F(z,w/A(z)^2,\Lambda_0)$ are holomorphic
functions of $z$ since this is the case for $H^j(z)$. This
implies that $H(z,w)=K(z,w,\Lambda_0)$, the holomorphic
mapping given in Step 4, and, hence, the proof of Theorem 4
is complete.\qed\enddemo

\heading \S7 Examples \endheading

In this section we shall present some examples which were
computed by implementing the algorithm described in the
previous section using  Mathematica for the symbolic
manipulation. All examples given here are hypersurfaces in
$\bC^2$. In Example 7.1,
$\Aut (M,p)$ is computed at a
$2$-nondegenerate (and hence Levi degenerate) point and the
group is not connected. Example 7.2 gives a family of
hypersurfaces
$M_{(a,b)}$, and we calculate $\scrF(M,0;M',0)$ for any pair
$M$, $M'$ in this family. The third example is a
hypersurface through 0 for which $\Aut(M,0)$ consists of
precisely two elements.  The
last example shows that even for an everywhere Levi
nondegenerate hypersurface, the real analytic set $\Sigma$
given in Corollary 4.1 need not be a manifold, although each
fiber $\Sigma_p$ must be, in view of Corollary 1.2.  Since
all the examples below are rigid hypersurfaces in $\bC^2$,
the calculation of $\hol_0 (M,0)$, the Lie algebra of $\Aut
(M,0)$ is essentially contained in Stanton [S1].  However,
the methods used in [S1] do not seem to lead to the
calculation of the entire stability group, even in these
examples.
\subhead Example 7.1
\endsubhead Let $M$ be the hypersurface given by
$$\im w =  (\re z) |z|^2.
$$ Note that the dilations $(z,w) \mapsto (tz,t^3w)$, with
$t\not=0$ real, are in $\Aut(M,0)$.  Calculations using the
algorithm show that these are the only elements of this
group. Hence the group has two connected components,
correstponding to $t>0$ and $t <0$.

\subhead Example 7.2
\endsubhead For $a\in \bR\setminus 0$ and $b\in
\bC\setminus 0$ let $M_{(a,b)}$ be the hypersurface given by
$$\im w = a|z|^2 + (\re b z) |z|^2.
$$ The biholomorphism $$(z,w)\mapsto \left(\frac a{\bar  b}
z,
\frac{a^3} {|b|^2} w\right)$$ maps $(M_{(1,1)},0)$ onto
$(M_{(a,b)},0)$. A calculation using the algorithm shows
that $\Aut (M_{(1,1)},0)$ consists of only the identity
map.  Hence
$\scrF(M_{(a,b)},0;M_{(a',b')},0)$ consists of a single
element for all $a,b,a',b'$ as above.

Let $L\subset \bC^2$ be the real line
$\{(0,x), x\in \bR\}$. A simple calculation shows that if
$M$ is the hypersurface in Example 7.1, then for every
$p\in M\setminus L$, there exists an
$\a > 0$ such that $(M,p)$ is biholomorphically equivalent
to $(M_{(\a,1)},0)$.  Combining Examples 7.1 and 7.2, we
find that for all $p,p' \in M\setminus L$, $(M,p)$ is
biholomorphically equivalent to $(M,p')$.  Furthermore,
$\Aut(M,p)$ consists of the dilations if $p \in L $ and
contains only the identity if $p \in M\setminus L$.
Observe that $M$ is 2-nondegenerate along $L$ and
$1$-nondegenerate elsewhere.  In the notation of Corollary
4.1 we may take
$k_0=2$ and conclude from the above that $\Sigma\subset \gp
{4} (\bC^2)$ consists of the transverse union of a three
dimensional manifold and a two dimensional one.

\subhead Example 7.3
\endsubhead Let $M$ be the hypersurface given by
$$\im w = |z|^2 + (\re z^2) |z|^2.
$$ Then $\Aut(M,0)$ consists of exactly two elements, namely
the identity and the map $(z,w) \mapsto (-z,w)$.

\subhead Example 7.4
\endsubhead Consider the hypersurface $M$ given by
$$\im w = |z|^2 +  |z|^4.
$$ The only elements of $\Aut(M,0)$ are the rotations
$(z,w) \mapsto (e^{i\theta}z,w)$ for
$\theta \in \bR$. A simple calculation shows that for
$a>0$ and $p_a = (a,a^2+a^4) \in M$ the germ $(M,p_a)$ is
biholomorphically equivalent to $(M_a,0)$, where
$M_a$ is defined by
$$\im w = (1+4a^2)|z|^2 +4\,a\,\re(z)|z|^2+  |z|^4.
$$ For $a = 1/2$, we were able to show that $\Aut(M_a,0)$
consists of only the identity. We write $\bold 1$ for the
identity in $\gp 2_0 (\bC^2)$. In the calculation of the
equations  for the set
$\Sigma_0 =
\eta(\Aut(M_{1/2},0))$  we could extract $d$ equations
defining $\Sigma_0$, $c_1(\L,\bar\L),
\ldots,c_d(\L,\bar\L)
$, where $d = \dim_\bR G^2_0(\bC^2)$, satisfying
$$ \text {\rm rank} \left(\nabla_{\L,\bar\L}c_1({\bold 1},
\bar {\bold 1}), \ldots,\nabla_{\L,\bar\L}c_d({\bold 1},
\bar {\bold 1})\right)=d.
$$ Since the defining equations vary analytically with the
parameter $a$, the rank of these equations is still $d$ for
$a$ near $1/2$.  From this we conclude that $\Sigma=
\sigma^2(\scrA(M))$, as defined in Corollary 4.1, is a
three dimensional manifold near the point
$(p_{1/2},{\bold 1})$ in the fiber bundle $G^2(\bC^2)$ over
$M$.  On the other hand, the set $\Sigma$ is not a three
dimensional manifold near the point $(0,{\bold 1})$, since
it contains the three dimensional manifold $(M,{\bold 1})$
and a circle over the point
$0\in M$.  Since $(p_{1/2}, {\bold 1})$  is connected to
$(0, {\bold 1})$ in $\Sigma$, the latter cannot be a
manifold.

The conclusion that $\Sigma$ is not a manifold could also
be deduced using Theorem 4, applied at the point
$p_{1/2}$ with the notation above, and some additional work.

\enddocument